\documentclass[11pt]{article}
\usepackage{amsmath}
\usepackage{amsfonts}
\usepackage{amssymb}
\usepackage{enumerate}
\usepackage[ruled,vlined,linesnumbered]{algorithm2e}
\usepackage{amsthm}
\usepackage{authblk}
\usepackage{graphicx}
\usepackage{verbatim,color,hyperref}
\newtheorem{myproposition}{Proposition}[section]
\newtheorem{mytheorem}[myproposition]{Theorem}
\newtheorem{mylemma}[myproposition]{Lemma}

\newtheorem{mycorollary}[myproposition]{Corollary}

\newtheorem{myproblem}[myproposition]{Problem}

\def\gr{\mathcal{G}}

\makeatletter
\def\imod#1{\allowbreak\mkern10mu({\operator@font mod}\,\,#1)}
\makeatother
\begin {document}

\title{Group Irregularity Strength of Connected Graphs}


\author[1,4]{Marcin Anholcer}
\author[2,4]{Sylwia Cichacz}
\author[3,4]{Martin Milani\u{c}}
\affil[1]{\scriptsize{}Pozna\'n University of Economics, Faculty of Informatics and Electronic Economy}
\affil[ ]{Al. Niepodleg{\l}o\'sci 10, 61-875 Pozna\'n, Poland, \textit{m.anholcer@ue.poznan.pl}}
\affil[ ]{}
\affil[2]{AGH University of Science and Technology, Faculty of Applied Mathematics}
\affil[ ]{Al. Mickiewicza 30, 30-059 Krak\'ow, Poland, \textit{cichacz@agh.edu.pl}}
\affil[ ]{}
\affil[3]{University of Primorska, UP IAM, Muzejski trg 2, SI6000 Koper, Slovenia}
\affil[ ]{martin.milanic@upr.si}
\affil[ ]{}
\affil[4]{University of Primorska, UP FAMNIT, Glagolja\v ska 8, SI6000 Koper, Slovenia}
\affil[ ]{\{marcin.anholcer, sylwia.cichacz-przenioslo, martin.milanic\}@famnit.upr.si}


\maketitle

\begin{abstract}
We investigate the \textit{group irregularity strength} ($s_g(G)$) of graphs, that is, the smallest value of $s$ such that taking any Abelian group $\gr$ of order $s$, there exists a function $f:E(G)\rightarrow \gr$ such that the sums of edge labels at every vertex are distinct. We prove that for any connected graph $G$ of order at least $3$, $s_g(G)=n$ if $n\neq 4k+2$ and $s_g(G)\leq n+1$ otherwise, except the case of some infinite family of stars.
\end{abstract}

\noindent\textbf{Keywords:} irregularity strength, graph labelling, Abelian group\\
\noindent\textbf{MSC:} 05C15, 05C78.

\section{Introduction}

It is a known fact that in any simple graph $G$ there are at least two vertices of the same degree. The situation changes if we consider multigraphs. Each multiple edge may be represented with some integer label and the (weighted) degree of any vertex $x$ is then calculated as the sum of labels over all edges incident to $x$. The maximum label $s$ is called the \textit{strength} of the labelling. The labelling itself is called \textit{irregular} if the weighted degrees of all the vertices are distinct. The smallest value of $s$ that allows some irregular labelling is called \textit{irregularity strength of $G$} and denoted  by $s(G)$.

The problem of finding $s(G)$ was introduced by Chartrand et al.~in \cite{ref_ChaJacLehOelRuiSab1} and investigated by numerous authors. Best published result due to Kalkowski et al.~(see \cite{ref_KalKarPfe1}) is $s(G)\leq 6n/\delta$. There are some signals that it was recently improved by Przyby{\l}o (\cite{ref_Prz3}) for dense graphs of sufficiently big order ($s(G)\leq 4n/\delta$ in this case). Exact value of $s(T)$ for a tree $T$ was investigated e.g. by Aigner and Triesch (\cite{ref_AigTri}), Amar and Togni (\cite{ref_AmaTog}), Ferrara et al.~(\cite{ref_FerGouKarPfe}) and Togni (\cite{ref_Tog1}).

On the other hand, numerous authors studied various labelling problems when elements of finite Abelian groups were used instead of integers to label either vertices or edges of graph. We give only few examples here. Graham and Sloane in \cite{ref_GraSlo} studied harmonious graphs, i.e., graphs for which there exists an injection $f:V(G)\rightarrow \mathbb{Z}_q$ that assigns to every edge $(x,y)\in E(G)$ unique sum $f(x)+f(y)$ modulo $q$. Beals et al.~(see \cite{ref_BeaGalHeaJun}) considered the concept of harmoniousness with respect to arbitrary Abelian groups. \.Zak in \cite{ref_Zak} generalized the problem and introduced new parameter, \textit{harmonious order of $G$}, the smallest number $t$ such that injection $f:V(G)\rightarrow \mathbb{Z}_t$ (or surjection if $t<V(G)$) produces distinct edge sums. Hovey in \cite{ref_Hov} considers the so-called $A-cordial$ labellings, where for a given Abelian group $A$ and a graph $G$ one wants to obtain such a vertex labelling that the classes of vertices labelled with one label are (almost) equinumerous and so are the classes of edges with the same sums. Cavenagh et al.~(\cite{ref_CavComNel}) consider \textit{edge-magic total labellings} with finite Abelian groups, i.e., the labelings of vertices and edges resulting in equal edge sums. Froncek in \cite{ref_Fro} defined the notion of group distance magic graphs, i.e., the graphs allowing the bijective labelling of vertices with elements of an Abelian group resulting in constant sums of neighbour labels. Stanley in \cite{ref_Sta} studied the vertex-magic labellings of edges with the elements of an Abelian group $A$, i.e., labellings, where the resulting weighted degrees are constant. Kaplan et al.~in \cite{ref_KapLevRod} considered vertex-antimagic edge labellings, i.e., the bijections $f:E(G)\rightarrow A$, where $A$ is a cyclic group, resulting in distinct weighted degrees of vertices.

The problem considered in this paper arises as the complement of the research conducted so far. Assume we are given an arbitrary graph $G$ of order $n$ with no components isomorphic to $K_1$ or $K_2$. Assume $\gr$ is an Abelian group  of order $m\geq n$ with the operation denoted by $+$ and neutral element $0$. For convenience we will write $ka$ to denote $a+a+\dots+a$ (where element $a$ appears $k$ times), $-a$ to denote the inverse of $a$ and we will use $a-b$ instead of $a+(-b)$. Moreover, the notation $\sum_{a\in S}{a}$ will be used as a short form for $a_1+a_2+a_3+\dots$, where $a_1, a_2, a_3, \dots$ are all the elements of the set $S$.

We define edge labelling $f:E(G)\rightarrow \gr$ leading us to the weighted degrees defined as the sums:
$$
w(v)=\sum_{e\ni v}{f(e)}
$$
We call $f$ \textit{$\gr$-irregular} if all the weighted degrees are distinct. The \textit{group irregularity strength} of $G$, denoted $s_g(G)$, is the smallest integer $s$ such that for every Abelian group $\gr$ of order $s$ there exists $\gr$-irregular labelling $f$ of $G$. The main result of our paper is the following theorem, determining the value of $s_g(G)$ for every connected graph $G$ of order $n\geq 3$.

\begin{mytheorem}\label{main_thm}
Let $G$ be arbitrary connected graph of order $n\geq 3$. Then
$$
s_g(G)=\begin{cases}
n+2&\text{when   } G\cong K_{1,3^{2q+1}-2} \text{  for some integer   }q\geq 1\\
n+1&\text{when   } n\equiv 2 \imod 4 \wedge G\not\cong K_{1,3^{2q+1}-2} \text{  for any integer   }q\geq 1\\
n&\text{otherwise}
\end{cases}
$$
\end{mytheorem}

We also show that the following theorem is true.

\begin{mytheorem}\label{main_thm2}
Let $G$ be arbitrary connected graph of order $n\geq 3$.
Then, for every $k>s_g(G)$ and every finite Abelian group $\gr$ of order $k$, $G$ admits a $\gr$-irregular labelling,
except for the cases when:
\begin{itemize}
\item
$G\cong K_{1,n-1}$ and $\gr\cong \mathbb{Z}_3\times \mathbb{Z}_3\times\dots\times \mathbb{Z}_3=(\mathbb{Z}_3)^q$ for some $q$ such that $3^q= n+1$
\item
$\gr\cong  \mathbb{Z}_2 \times  \mathbb{Z}_2\times \ldots \times \mathbb{Z}_2=(\mathbb{Z}_2)^q$ for some $q$ such that $2^q=n+2$.
\end{itemize}
\end{mytheorem}

\section{Proof of Theorem \ref{main_thm}}\label{sec:main}

In order to distinguish $n$ vertices in the graph we need at least $n$ distinct elements of $\gr$. However, $n$ elements are not always enough, as shows the following lemma.

\begin{mylemma}\label{lemma_below}
Let $G$ be of order $n$, if $n \equiv 2 \imod 4$, then $s_g(G) \geq n+1$.
\end{mylemma}

\noindent\textbf{Proof.}
Let $\gr$ be an Abelian group of order $n=2(1+2k)$. The fundamental theorem of finite abelian groups states that the finite abelian group $\gr$ can be expressed as the direct sum of cyclic subgroups of prime-power order. This implies that $\gr \cong \mathbb{Z}_2 \times \mathbb{Z}_{p_{1}^{\alpha_1}}\times \mathbb{Z}_{p_{2}^{\alpha_2}}\times \ldots \times \mathbb{Z}_{p_{m}^{\alpha_m}}$,
where $n=2\prod_{i=1}^m{p_i^{\alpha_i}}$ and $p_i>2$ for $i=1,\dots,m$ are not necessarily distinct primes. This implies that one can write $a \in \gr$ as $a=(a_0,a_1,\ldots,a_m)$. Notice that in the group $\gr$ we have $1+2k$ elements with the first coordinate $0$ and
$1+2k$ with the first coordinate $1$.
Let now $w(G)=\sum_{v \in V(G)}w(v)=\sum_{a\in \gr}a$. Observe that $w(G)$ is a vector with the first coordinate $1$ (since we are summing in $\mathbb{Z}_2$).
On the other hand $w(G)=\sum_{v \in V(G)}(\sum_{ev}f(e))$, so each label $f(e)$ for any $e \in E(G)$ appears in the sum twice. Therefore
 $w(G)$ is a vector with the first coordinate $0$ (since we are summing in $\mathbb{Z}_2$), contradiction.

\qed

We continue with the following lemma, determining the group irregularity strength of stars.

\begin{mylemma}\label{lemma_stars}
Let $K_{1,n-1}$ be a star with $n-1$ pendant vertices and  $n\geq 3$. Then
$$
s_g(K_{1,n-1})=\begin{cases}
n+2&\text{when   }n\equiv 2 \imod 4 \wedge n=3^q-1\text{  for some integer   }q\geq 1\\
n+1&\text{when   }n\equiv 2 \imod 4 \wedge n\neq 3^q-1\text{  for any integer   }q\geq 1\\
n&\text{otherwise}.
\end{cases}
$$
\end{mylemma}

\noindent\textbf{Proof.} If $n$ is odd then we put all the elements of $\gr$ other than $0$ on the pendant edges and obtain this way distinct weighted degrees (same as edge labels) on the leafs and weighted degree $0$ in the central vertex.

If $n=4k$ for some $k\geq 1$, then we distinguish two cases, depending on the number of involutions. If there is only one involution, then it is  guaranteed that there exists a subgroup of $\gr$ isomorphic with $\mathbb{Z}_4$: $\{0,a,2a,3a\}$ ($2a$ is the only involution here). In such a situation we label the edges with all the elements of $\gr$ except $3a$, assigning this way the same values to the weighted degrees of all the pendant vertices. It is straightforward to check that the weighted degree of central vertex is $3a$. If there are at least two involutions, then the sum of all the elements of $\gr$ is $0$ (see e.g. \cite{ref_ComNelPal}, Lemma 8). We put on the edges all the elements of $\gr$ but $0$ and thus obtain distinct weighted degrees of pendant vertices not equal to $0$ and weighted degree $0$ of the central vertex.

If $n=4k+2$ for some natural $k\geq 1$, the order of $\gr$ must be at least $4k+3$ by Lemma~\ref{lemma_below}. Assume that $|\gr|=4k+3$, then there is no involution in $\gr$. If there is an element $a$ in $\gr$ of order more than $3$, then we assign to three edges labels $a$, $-2a$ and $0$ and we put $2k-1$ pairs $\{a_j,-a_j\}$, where $a_j\not\in\{0,a,-a,2a,-2a\}$, on the remaining edges, obtaining this way the $\gr$-irregular labelling.

If all the elements of $\gr$ have order $3$, then $n=3^q-1$ and $\gr\cong \mathbb{Z}_3\times \mathbb{Z}_3\times\dots\times \mathbb{Z}_3$. Assume we are able to label $K_{1,n-1}$ with $n+1$ labels from $\gr$. In such a situation we would have to use $n-1$ distinct elements of $\gr$ on the edges, which would leave us two distinct elements, say $a_1$ and $a_2$. The weighted degree of the central vertex would be $-(a_1+a_2)$. This number should be distinct from all other degrees, so one of the equalities $-(a_1+a_2)=a_1$ or $-(a_1+a_2)=a_2$ should be satisfied. In both cases it follows that $a_1=a_2$, contradiction. Thus we have to use the group $\gr$ of order at least $n+2$ in order to obtain $\gr$-irregular labelling. Observe that $n+2$ can not be equal to $2^p$ for any natural number $p$ (this follows from the Mih\v{a}ilescu Theorem, also known as the Catalan Conjecture, see \cite{ref_Mih}). Thus in any group of order $n+2$ there are two distinct elements $a_1$ and $a_2$ such that $a_1+a_2=0$. If there is more than one involution in $\gr$, then we label the edges with all the elements of $\gr$ but $0$, $a_1$ and $a_2$ and obtain this way the sum $0$ at the central vertex, distinct from all the other weighted degrees. If there is exactly one involution $i$ in $\gr$, then $\gr$ has a subgroup $\gr_1$ isomorphic with $\mathbb{Z}_4$: $\gr_1=\{0,a,2a=i,3a\}$. Thus we can put $0$, $a$ and $2a$ on three of the edges of the star. The remaining $4k$ elements of $\gr$ form $2k$ distinct pairs $\{a_j,-a_j\}$ such that $a_j\not\in\gr_1$, so we can put $2k-1$ of them on the remaining edges. Finally the central vertex obtains the weight $3a$ and all the vertex weights are distinct.

\qed

Now we are going to determine the value of $s_g(K_{1,n-1})$ for arbitrary tree $T$ not being a star.

\begin{mylemma}\label{lemma_trees}
Let $T$ be arbitrary tree on $n\geq 3$ vertices not being a star. Then
$$
s_g(T)=\begin{cases}
n+1&\text{when   }n\equiv 2 \imod 4 \\
n&\text{otherwise}.
\end{cases}
$$
\end{mylemma}

\noindent\textbf{Proof.} Assume we coloured properly the vertices of $T$ with two colours, obtaining colour classes $V_1$ and $V_2$. Given any two vertices $x_1$ and $x_2$, there exists unique path $P(x_1,x_2)$ that joins them. If $x_1$ and $x_2$ belong to the same colour class, then the path consists of odd number of vertices and even number of edges (we will call such a path \textit{odd path}). If $x_1$ and $x_2$ belong to distinct colour classes, then the path consists of even number of vertices and odd number of edges (we will call such a path \textit{even path}).

We start with $0$ on all the edges of $T$. Then, in every step, we will add some labels to all the edges of chosen path $P(x_1,x_2)$. To be more specific, we will add some label $a$ to all the edges having odd position on the path (starting from $x_1$) and $-a$ to all the edges having even position. We will denote such situation by $\phi(x_1,x_2)=a$. Observe that if $P(x_1,x_2)$ is odd, then putting $\phi(x_1,x_2)=a$ increases the weighted degree of $x_1$ by $a$ and the weighted degree of $x_2$ by $-a$. If $P(x_1,x_2)$ is even, then the weighted degrees of $x_1$ and $x_2$ increase by $a$. In both cases the weighted degrees of the remaining vertices stay unchanged.

Let us start with the case when $n=2k+1$ for some integer $k\geq 1$. Let $\gr$ be an arbitrary Abelian group of order $n$. As $\gr$ does not have any elements of order $2$, we can choose $k$ elements $a_1, a_2,\dots, a_k\in \gr$ such that $a_i\not\in\{a_j,-a_j\}$ for $i\neq j$ and $a_i\neq 0$ for $1\leq i\leq k$. One of the colour classes of $V(T)$, say $V_1$, has odd number of vertices and $V_2$  even. We join the vertices of $V_2$ in pairs, then do the same with all the vertices in $V_1$ but one, say $x_0$. We obtained this way exactly $k$ monochromatic pairs $(x_{j,1},x_{j,2})$. Now we put $\phi(x_{j,1},x_{j,2})=a_j$, for $1\leq j\leq k$. This way we obtain the $\gr$-irregular weighting: $w(x_0)=0$, $w(x_{j,1})=a_j$ for $1\leq j\leq k$, $w(x_{j,2})=-a_j$ for $1\leq j\leq k$.

If $n=4k$ for some integer $k\geq 1$, then we distinguish two cases. If there is only one involution in $\gr$, then there exists a subgroup $\{0, a, 2a, 3a\}$ of $\gr$, where $2a$ is the only element of order $2$ in $\gr$. In such a case we choose two vertices $x_1$, $x_2$ from one of the color classes and one vertex $x_0$ from the other one. We put $\phi(x_0,x_1)=a$ and $\phi(x_0,x_2)=2a$, obtaining this way $w(x_1)=a$, $w(x_2)=2a$ and $w(x_0)=3a$. The number of the remaining vertices in one of the color classes is now odd and in the second one even. Thus we can proceed as in the case of $n$ odd, using the remaining labels $\{a_j,-a_j\}$ such that $a_j\not\in\{a,2a,3a\}$ and obtaining this way $\gr$-irregular labelling of $T$. If there are more involutions $a_1, a_2,\dots, a_r$, then their number $r$ is odd and their sum equals $0$ (see e.g. \cite{ref_ComNelPal}, Lemma 8). If $r\leq n/2$, then we choose one vertex $x_0$ from the colour class with less vertices (say $V_p$) and $r$ vertices $x_1$, $x_2$, $\dots$, $x_r$ from $V_{3-p}$ and we put $\phi(x_0,x_j)=a_j$ for $j=1,\dots,r$. This way we obtain $w(x_0)=0$ and $w(x_j)=a_j$ for $j=1,\dots,r$. If the numbers of vertices in $V_1$ and $V_2$ are both even, we continue like in the case of $n$ odd (this time we do not obtain $w(x)=0$ for any new vertex). If both numbers are odd, we choose one vertex $x_{r+1}$ from $V_p$ and any element $a_{r+1}\not\in\{0,a_1,\dots, a_r\}$. By putting $\phi(x_0,x_{r+1})=a_{r+1}$ we obtain finally $w(x_0)=a_{r+1}$ and $w(x_{r+1})=-a_{r+1}$. Now the number of remaining vertices in $V_p$ is even and in $V_{3-p}$ odd. As we finally did not assign weighted degree $0$ to any vertex, we can proceed as in the case of $n$ odd. Last case to analyse here is when $j>n/2$. But in such situation $\gr\cong \mathbb{Z}_2\times \mathbb{Z}_2\times \dots\times \mathbb{Z}_2$ and all the elements of $\gr$ but $0$ have order $2$. In such case we chose any vertex $x_0$ of $T$ and put $\phi(x_0,x_j)=a_j$ for $j=1,\dots, n-1$ for distinct elements $a_j\neq 0$. This way we obtain $w(x_j)=a_j\neq 0$ for $j=1,\dots, n-1$ and $w(x_0)=0$.

If $n=4k+2$, then we will use the group of order $n+1$. We have to distinguish two cases. If both colour classes of $T$ are even, then we proceed as in the case when $n$ is odd (the difference is that there is no vertex $x_0$ with $w(x_0)=0$). If both colour classes are odd, we have to start with reducing their sizes in such a way that they both become even.

If there is an element of $\gr$ of order greater than $3$, say $a$, then we select three arbitrary vertices $x_1$, $x_2$, $x_3$ from one colour class and any vertex $x_0$ from the other one and we put $\phi(x_1,x_0)=a$, $\phi(x_2,x_0)=-2a$, $\phi(x_3,x_0)=0$. This way we obtain the following weighted degrees: $w(x_0)=-a$, $w(x_1)=a$, $w(x_2)=-2a$, $w(x_3)=0$. As we can easily see, these degrees are distinct and we still have $k-2$ pairs of elements $\{a_j,-a_j\}$ to label the remaining vertices as in the previous cases.

If all the elements of $\gr$ have order $3$, then $n\geq 26$ (in fact, we need only $n\geq 10$). We choose $a, b, c\in \gr$ such that $a\neq 0$, $b\neq 0$, $c\neq 0$, $a\not\in\{b,-b\}$, $c\not\in\{a, -a, b, -b, a+b, -(a+b), a-b, b-a\}$. As $T$ is not star, we can choose five vertices $x_1$, $x_2$, $x_3$, $x_4$, $x_5$ from one colour class and three $y_1$, $y_2$, $y_3$ from another one. Now we put $\phi(x_1, y_1)=a$, $\phi(x_2, y_1)=a$, $\phi(x_2, y_2)=2a+b$, $\phi(x_3, y_2)=a+b$, $\phi(x_3, y_3)=2a+2b+c$, $\phi(x_4, y_3)=a+b+c$, $\phi(x_5, y_3)=0$. This way we obtain $8$ distinct weighted degrees $w(x_1)=a$, $w(x_2)=b$, $w(x_3)=c$, $w(x_4)=a+b+c\not\in\{0,a,b,c,-a,-b,-c\}$, $w(x_5)=0$, $w(y_1)=-a$, $w(y_2)=-b$, $w(y_3)=-c$ and we still have $2k-3$ pairs $\{a_j,-a_j\}$ to distinguish remaining even numbers of vertices in both colour classes.

\qed

The main result of our paper follows easily from the above lemmas. If $G$ is a star, then we use lemma \ref{lemma_stars}. Otherwise we choose any spanning tree of $G$ not being a star, and use lemmas \ref{lemma_trees} and \ref{lemma_below}, labelling all the remaining edges with $0$.

\section{Proof of Theorem \ref{main_thm2}}\label{sec:main2}

Before we prove Theorem \ref{main_thm2}, we need  the following technical lemma that was proved in~\cite{ref_AnhCic2}.

\begin{mylemma}[\cite{ref_AnhCic2}]\label{lemma_inv}
Let $\gr$ be an Abelian group with involutions set $I^\star=\{i_1,\dots, i_{2^k-1}\}$, $k\geq 2$ and let $I=I^\star\cup \{0\}$. Then for any given $r$ such that $0\leq r \leq 2^k$, there exists set $R\subseteq I$,  $|R|=r$, such that
$$
\sum_{i\in R}{i}=0
$$
if and only if $r\not\in \{2,2^k-2\}$.
\end{mylemma}

In the next step we are going to show the construction of desired labelling for stars.

\begin{mylemma}\label{duzes}
Let $K_{1,n-1}$ be a star with $n-1$ pendant vertices and  $n\geq 3$. Then  $K_{1,n-1}$ admits $\gr$-irregular labelling for any finite Abelian group $\gr$ of order $k> s_g(K_{1,n-1})$ except the cases when $\gr\cong \mathbb{Z}_3\times \mathbb{Z}_3\times\dots\times \mathbb{Z}_3=(\mathbb{Z}_3)^q$ for some $q$ such that $3^q= n+1$ and $\gr\cong  \mathbb{Z}_2 \times  \mathbb{Z}_2\times \ldots \times \mathbb{Z}_2=(\mathbb{Z}_2)^q$ for some $q$ such that $2^q=n+2$.
\end{mylemma}
\noindent\textbf{Proof.} We can write $k=2^{p}(2m+1)$ for some natural numbers $p$ and $m$.

Suppose first that there exists at most one involution $i \in \gr$ or $2m \geq n-1$. If $n$ is odd then we put $\frac{n-1}{2}$ pairs $\{a_j,-a_j\}$  ($a_j \neq 0$) on the pendant edges and obtain this way distinct weighted degrees (same as edge labels) on the leafs and the weighted degree $0$ in the central vertex.

If $n$ is even and there is an element $a\in \gr$ of order more than $3$, then  we assign to three edges labels $a$, $-2a$ and $0$ and we put $\frac{n-4}{2}$ pairs $\{a_j,-a_j\}$, where $a_j\not\in\{0,a,-a,2a,-2a\}$, on the remaining edges, obtaining this way the $\gr$-irregular labelling of $G$ (such number of pairs exist, as $k>n$).

If $n$ is even and all the elements of $\gr$ have order less than $4$ and there exists the involution $i \in \gr$ then $k\geq n+2$ (as $k$ is even) and we assign to three edges labels $a$, $i$ and $0$ and we put $\frac{n-4}{2}$ pairs $\{a_j,-a_j\}$, where $a_j\not\in\{0,a,-a,i,a+i,-a+i\}$, on the remaining edges, obtaining this way the $\gr$-irregular labelling.

Finally, if all the elements of $\gr$ have order $3$, then $k=3^r$ for some $r$ and $\gr\cong \mathbb{Z}_3\times \mathbb{Z}_3\times\dots\times \mathbb{Z}_3$. If now  $n+1=k$, then $s_g(G)=n$ and $n\equiv 0 \imod 4$ by Theorem~\ref{main_thm}. Assume that we are able to label $K_{1,n-1}$ with $n+1$ labels from $\gr$. In such a situation we would have to use $n-1$ distinct elements of $\gr$ on the edges, which would leave us two distinct elements, say $a_1$ and $a_2$. The weighted degree of the central vertex would be $-(a_1+a_2)$. This number should be distinct from all other degrees, so one of the equalities $-(a_1+a_2)=a_1$ or $-(a_1+a_2)=a_2$ should be satisfied. In both cases it follows that $a_1=a_2$, contradiction. Thus we have to use the group $\gr$ of order at least $n+2$ in order to obtain $\gr$-irregular labelling. Observe that $n+2$ can not be equal to $2^p$ for any natural number $p$ (this follows from the Mih\v{a}ilescu Theorem, also known as the Catalan Conjecture, see \cite{ref_Mih}). Thus in any group of order $n+2$ there are two distinct elements $a_1$ and $a_2$ such that $a_1+a_2=0$. We label the edges with all the elements of $\gr$ but $0$, $a_1$ and $a_2$ and obtain this way the sum $0$ at the central vertex, distinct from all the other weighted degrees. Thus we can assume that $n+3\leq k$ and we assign to three edges labels $a$, $b\neq 2a$ and $2a+2b$ and we put $\frac{n-4}{2}$ pairs $\{a_j,-a_j\}$, where $a_j\not\in\{0,a,2a,b,2b,a+b,2a+2b\}$, on the remaining edges, obtaining this way the $\gr$-irregular labelling.

Suppose now that there exist $t>1$ involutions $i_1,\ldots,i_t$ in $\gr$. Recall that $t=2^q-1 $ for some $1<q\leq p$ and $\sum_{j=1}^ti_t=0$ (see e.g. \cite{ref_ComNelPal}, Lemma 8). Let $I^\star$ denote the set of all the involutions and let $I=I^\star\cup\{0\}$.

If $t\leq n $, then in the case $n$ being even (odd) we put $i_j$ on $j=1,\dots,t$ (respectively $j=1,\dots,t-1$) on $t$ (respectively $t-1$) edges and $(n-t-1)/2$ ($(n-t)/2$, respectively) distinct pairs $\{x_l,-x_l\}$ on the remaining edges.  We obtain this way distinct weighted degrees (same as edge labels) on the leafs and the weighted degree $0$ ($i_t$, respectively) in the central vertex.

Assume now that $t=n+1$, therefore $n$ is even. If there exists an element $a\in \gr$ such that $2a \neq 0$, then we assign to two edges labels $a$, $-a$. Using Lemma~\ref{lemma_inv}, we can choose $n-2=2^{q}-4$ elements $i_{j_1},i_{j_2},\ldots,i_{j_{n-2}}\in I$, such that

$$\sum_{l=1}^{n-2}i_{j_l}=0.$$

We put the elements $i_{j_1},i_{j_2},\ldots,i_{j_{n-3}}$ on the remaining edges, obtaining this way the $\gr$-irregular labelling. If there is no such element $a \in \gr$, then $\gr\cong  \mathbb{Z}_2 \times  \mathbb{Z}_2\times \ldots \times \mathbb{Z}_2$ and $k=2^p=n+2$. Let us assume that we managed to distinguish all the vertices and we did not use labels $a,b \in \gr$ ($a\neq b$). Thus $\sum_{v\in V(G)}w(v)=\sum_{g\in \gr}g-a-b=-a-b \neq 0$. On the other hand each label $f(e)$ for any $e \in E(G)$ appears in the sum twice. Therefore  $\sum_{v\in V(G)}w(v)= 0$. The contradiction shows that it is impossible to find desired labelling if $\gr\cong  \mathbb{Z}_2 \times  \mathbb{Z}_2\times \ldots \times \mathbb{Z}_2$ and $2^p=n+2$.

Let us consider now the case $t\geq n+2$. We have that $2<n<  2^p-2$, thus using Lemma~\ref{lemma_inv} we can choose $n$ elements $i_{j_1},i_{j_2},\ldots,i_{j_{n}}\in I$ such that

$$\sum_{l=1}^{n}i_{j_l}=0.$$

We put the elements $i_{j_1},i_{j_2},\ldots,i_{j_{n-1}}$  on the edges obtaining this way distinct weighted degrees (same as edge labels) for the leafs and the weighted degree $i_n$ for the central vertex. \qed

\begin{mylemma}\label{duzed}
Let $T$ be arbitrary tree on $n\geq 4$ vertices not being a star. Then  $T$ admits $\gr$-irregular labelling for any abelian group $\gr$ of order $k> s_g(T)$ except the case when $\gr\cong  \mathbb{Z}_2 \times  \mathbb{Z}_2\times \ldots \times \mathbb{Z}_2=(\mathbb{Z}_2)^q$ for some $q$ such that $2^q=n+2$.
\end{mylemma}

\noindent\textbf{Proof.}
We can  write $k=2^{p}(2m+1)$ for some natural numbers $p$ and $m$.

Suppose first that there exists exactly one involution $i \in \gr$ or $2m \geq n-1$.
If now  $n$ is odd or both colour classes of $T$ are even, then since we have at least $\lfloor \frac{n}{2}\rfloor$ pairs $\{a_i,-a_i\}$ such that $2a_i \neq 0$, we join the vertices of $V_1$ in pairs (except possibly one vertex if $|V_1|$ is odd), then do the same with all the vertices in $V_2$. We obtain this way exactly $\lfloor\frac{n}{2}\rfloor$ monochromatic pairs $(x_{j,1},x_{j,2})$ plus possibly one additional vertex, say $x_0$. Now we put $\phi(x_{j,1},x_{j,2})=a_j$, for $1\leq j\leq \lfloor \frac{n}{2}\rfloor$. This way we obtain the $\gr$-irregular weighting, as $w(x_{j,1})=a_j=-w(x_{j,2})$ for $1\leq j\leq k$ and if $n$ is odd, then $w(x_0)=0$.

Consider now the case when both color classes are odd.

If there is an element $a\in\gr$ of order greater than $3$ or there exists involution $i \in \gr$ (thus   $k\geq n+2$), then we select three arbitrary vertices $x_1$, $x_2$, $x_3$ from one colour class and any vertex $x_0$ from the other one. In the first situation we put $\phi(x_1,x_0)=a$, $\phi(x_2,x_0)=-2a$, $\phi(x_3,x_0)=0$. This way we obtain the following weighted degrees: $w(x_0)=-a$, $w(x_1)=a$, $w(x_2)=-2a$, $w(x_3)=0$. In the second situation we choose some $a\not \in \{0,i\}$ and we put $\phi(x_1,x_0)=i$, $\phi(x_2,x_0)=a$, $\phi(x_3,x_0)=0$. This way we obtain the following weighted degrees: $w(x_0)=a+i$, $w(x_1)=i$, $w(x_2)=a$, $w(x_3)=0$. As we can easily see, these degrees are distinct and since $n$ is even and $k\geq n+1$ ($k\geq n+2$, respectively) we still have at least $\frac{n-4}{2}$ pairs of elements $\{a_j,-a_j\}$ to label the remaining vertices as in the previous cases.

If all the elements of $\gr$ have order $3$ and $k\geq 10$, then  we choose $a, b, c\in \gr$ such that $a\neq 0$, $b\neq 0$, $c\neq 0$, $a\not\in\{b,-b\}$, $c\not\in\{a, -a, b, -b, a+b, -(a+b), a-b, b-a\}$. As $T$ is not star, we can choose five vertices $x_1$, $x_2$, $x_3$, $x_4$, $x_5$ from one colour class and three $y_1$, $y_2$, $y_3$ from the other one. Now we put $\phi(x_1, y_1)=a$, $\phi(x_2, y_1)=a$, $\phi(x_2, y_2)=2a+b$, $\phi(x_3, y_2)=a+b$, $\phi(x_3, y_3)=2a+2b+c$, $\phi(x_4, y_3)=a+b+c$, $\phi(x_5, y_3)=0$. This way we obtain $8$ distinct weighted degrees $w(x_1)=a$, $w(x_2)=b$, $w(x_3)=c$, $w(x_4)=a+b+c\not\in\{0,a,b,c,-a,-b,-c\}$, $w(x_5)=0$, $w(y_1)=-a$, $w(y_2)=-b$, $w(y_3)=-c$ and we still have $\frac{k-9}{2}\geq\lceil\frac{n-8}{2}\rceil$ pairs $\{a_j,-a_j\}$ to distinguish remaining even numbers of vertices in both colour classes. If $n = 8$  and $\gr \cong \mathbb{Z}_3 \times \mathbb{Z}_3$, then since $T$ is not star, we can choose five vertices $x_1$, $x_2$, $x_3$, $x_4$, $x_5$ from one colour class and three $y_1$, $y_2$, $y_3$ from another one. Now we put $\phi(x_1, y_1)=(1,0)$, $\phi(x_2, y_1)=(2,0)$, $\phi(x_2, y_2)=(0,0)$, $\phi(x_3, y_2)=(1,1)$, $\phi(x_3, y_3)=(2,1)$, $\phi(x_4, y_3)=(2,1)$, $\phi(x_5, y_3)=(2,2)$. This way we obtain $8$ distinct weighted degrees $w(x_1)=(1,0)$, $w(x_2)=(2,0)$, $w(x_3)=(0,2)$, $w(x_4)=(2,1)$, $w(x_5)=(2,2)$, $w(y_1)=(0,0)$, $w(y_2)=(1,1)$, $w(y_3)=(0,1)$.

Suppose now that there exist $t>1$ involutions in $\gr$. Recall that $t=2^q-1 $ for some $1<q\leq p$. If $t\leq n$, then we choose $t$ vertices $x_1,x_2,\ldots, x_t\in V_1\cup V_2$ in such a way, that at least one of the numbers of remaining vertices in $V_1$ and $V_2$ is even. We put $\phi(x_1,x_j)=i_j$ for $j=2,\dots,t$ obtaining $w(x_j)=i_j$ for $j=1,\dots,t$. Since the numbers of remaining vertices in at least one of the colour classes $V_1$ and $V_2$ are even, we construct $\lfloor\frac{n-t}{2}\rfloor$ monochromatic pairs and use the pairs $\{x_l,-x_l\}$ to label them. If there is some unpaired vertex, then its weighted degree is $0$, so the obtained labelling is $\gr$-irregular.

Assume now that $t=n+1$. If there exists an element $a\in \gr$ such that $2a \neq 0$, then we  choose two vertices $x_{n-1}$, $x_{n}$ from one colour class and we put $\phi(x_{n-1}, x_{n})=a$. Using Lemma~\ref{lemma_inv} we can choose $t-3=n-2=2^{q}-4$ elements $i_{j_1},i_{j_2},\ldots,i_{j_{n-2}}\in I$ such that

$$\sum_{l=1}^{n-2}i_{j_l}=0.$$

We put $\phi(x_1,x_l)=i_{j_l}$ for $l=2,\dots, n-2$. This way we obtain $w(x_l)=i_{j_l}$ for $l=1,\dots,n-2$. If there does is no such element $a \in \gr$, then $\gr\cong  \mathbb{Z}_2 \times  \mathbb{Z}_2\times \ldots \times \mathbb{Z}_2$ and $k=2^p=n+2$. As in the proof of Lemma~\ref{duzes} let  us assume that we distinguished all vertices and we did not use labels $a,b \in \gr$ ($a\neq b$). Thus $\sum_{v\in V(G)}w(v)=\sum_{g\in \gr}g-a-b=-a-b \neq 0$. On the other hand each label $f(e)$ for any $e \in E(G)$ appears in the sum twice. Therefore  $\sum_{v\in V(G)}w(v)= 0$. This contradiction shows that it is impossible to find desired labelling in such a case.

Let us consider now the case $t\geq n+2$. We have $2<n<  2^q-2$, thus using Lemma~\ref{lemma_inv} we can choose $n$ elements $i_{j_1},i_{j_2},\ldots,i_{j_{n}}$ such that

$$\sum_{l=1}^{n-1}i_{j_l}=0.$$

We put $\phi(x_1,x_l)=i_{j_l}$ for $l=2,\dots, n$. This way we obtain $w(x_l)=i_{j_l}$ for $l=1,\dots,n$. \qed

Theorem~\ref{main_thm2} follows easily from the above lemmas. If $G$ is a star, then we use Lemma \ref{duzes}. Otherwise we choose any spanning tree of $G$ not being a star, and use Lemma \ref{duzed}  labelling all the remaining edges with $0$. Observe that in the latter case same argument as in the proof of Lemma \ref{duzed} shows that for every graph $G$ it is impossible to find $\gr$-irregular labelling of $G$ if $\gr\cong(\mathbb{Z}_2)^q$ for some $q$ such that $2^q=n+2$.

\section{Computational Complexity Issues}

In Section~\ref{sec:main}, the group irregularity strength $s_g(G)$ of an arbitrary connected graph $G$ was determined.
It seems natural to ask about the computational complexity of the corresponding  problem, where we assume that the group $\gr$ is given in the most compact form that follows from the fundamental theorem of Abelian groups.\\

\bigskip
\noindent {\sc Irregular Labeling}.\\
{\bf Instance:}
A connected graph $G$, an Abelian group $\gr$ of order $|\gr|$, given as the list of
(orders of) cyclic subgroups $\gr_1,\ldots, \gr_k$ of prime-power order such that $\gr\cong\gr_1\times\ldots\times \gr_k$.\\
{\bf Task:} Find a $\gr$-irregular labelling of $G$ or answer that it is impossible.

\bigskip
The proofs in Section~\ref{sec:main} and Section~\ref{sec:main2} are constructive and lead to efficient algorithms for the {\sc Irregular Labeling} problem.

In order to construct an irregular labeling whose existence is proved in Section~\ref{sec:main}, we first construct, in time $O(|V(G)|+|E(G)|)$,
a spanning tree $T$ of the given graph $G$. A proper $2$-coloring of $T$ can also be obtained in time $O(|V(T)|)= O(|V(G)|)$.
The rest of the construction of the irregular labeling reduces to solving constantly many instances of the following problem:
{\em Given a tree $T= (V,E)$ and a non-empty set $A\subseteq V$ of even cardinality, partition the elements of $A$ into pairs, say
$\{a_1,a_2\},\ldots, \{a_{2r-1},a_{2r}\}$, and construct the corresponding paths
 $P(a_1,a_2),\ldots, P(a_{2r-1},a_{2r})$ joining them.}
The problem is clearly solvable in polynomial time. However, if the set $A$ is large (as is the case for the subproblems one needs to solve in order to
construct an irregular labeling), the total length of the paths $P(a_1,a_2),\ldots, P(a_{2r-1},a_{2r})$ can be of the order $\Omega(|V|^2)$. Can one do better? In particular, can the {\sc Irregular Labeling} problem be solved in linear time?
We will show in this section that this is indeed the case.
With this goal in mind, we introduce the following optimization problem:

\bigskip
\noindent {\sc Shortest Path Collection}.\\
{\bf Instance:} A tree $T= (V,E)$ and a non-empty set $A\subseteq V$ of even cardinality.\\
{\bf Task:} Find a partition of the elements of $A$ into pairs, say
$\{\{a_1,a_2\},\ldots, \{a_{2r-1},a_{2r}\}\}$
such that  the value of $$\sum_{i = 1}^r{\rm dist}_T(a_{2i-1},a_{i})$$ is minimized.

\medskip
Here, the distance ${\rm dist}_T(\cdot, \cdot)$ is the usual graph-theoretic distance between vertices, that is, the number of edges on a shortest path connecting the two vertices.  The {\sc Shortest Path Collection} problem can be solved in polynomial time not only for trees but also for general graphs.
In fact, it can be solved in time $O(|V||E|+r^4)$ where $|A| = 2r$, by first computing in time $O(|V||E|)$
all pairwise vertex distances (this can be done using breadth-first search), and then solving an instance of the minimum weight perfect matching problem in a complete graph with vertex set $A$ and edge weights given by $w(xy) = {\rm dist}_T(x,y)$ for all pairs of distinct vertices $x,y\in A$ (this can be done, e.g., using the algorithm by Edmonds~\cite{ref_Edm}). As we show below, the problem can be solved in linear time for trees.

\begin{myproposition}
There exists a linear time algorithm for the {\sc Shortest Path Collection} problem.
Moreover, an optimal collection ${\cal P}$  of shortest paths, each connecting one pair of vertices from the partition of $A$,
can be constructed in linear time.
\end{myproposition}

\begin{proof}
Let $(T,A)$ be an instance to the {\sc Shortest Path Collection} problem.
The problem can be solved using a greedy algorithm, traversing the given tree $T=(V,E)$ bottom up and constructing optimal paths along the way.
The algorithm outputs  a pair $({\cal A}, {\cal P})$, where ${\cal A}$ is a partition of the elements of $A$ into pairs $\{a_1,a_2\},\ldots, \{a_{2r-1},a_{2r}\}$, and ${\cal P}$ is the collection of corresponding paths in $T$ connecting the paired vertices.
Traversing the tree bottom up, the algorithm updates two collections ${\cal P}$ and ${\cal Q}$
of paths in $T$ such that:
\begin{enumerate}[(1)]
  \item every path from ${\cal P}$ has exactly two endpoints in $A$ (in particular, ${\cal P}$ does not contain any trivial, one-vertex paths),
  \item every path from ${\cal Q}$ has exactly one endpoint in $A$ (in particular, ${\cal Q}$ can  contain several trivial  paths),
  \item no path from ${\cal Q}$ has a vertex in common with another path in ${\cal P}\cup {\cal Q}$, and
  \item every two paths in ${\cal P}$ are edge disjoint, and have at most one vertex in common.
\end{enumerate}
The set ${\cal Q}$ contains all paths that will be eventually extended to a path in the final solution ${\cal P}$.
Moreover, for each path $Q\in {\cal Q}$, its endpoints are denoted by $a(Q)$, $b(Q)$ in such a way that $a(Q)\in A$.
(If $Q$ is a trivial, one-vertex path then $a(Q) = b(Q)$.)

At every step of the algorithm, a vertex, say $v_i$, of $T$ is visited. Paths of ${\cal Q}$
that have a child of $v_i$ as one of their endpoints are greedily paired and merged, using vertex $v_i$, to form paths in ${\cal P}$.
At the end of this pairing procedure, one path from ${\cal Q}$
that has a child of $v_i$ as one of its endpoints could be left unpaired, in which case we extend it by the edge
connecting one of its endpoints to $v_i$. If $v_i$ belongs to $A$, this extended path is moved from ${\cal Q}$ to ${\cal P}$.
On the other hand, if all paths have been paired, then we check whether $v_i$ belongs to $A$ and if this is the case, we add to ${\cal Q}$ the trivial one-vertex path consisting of $v_i$.

Every time a path, say $P$, is added to the set ${\cal P}$, the set ${\cal A}$ is augmented with the pair containing the two endpoints of $P$.
Moreover, the algorithm keeps at every vertex a Boolean variable $q(v)$ such that $q(v) = 1$ if and only if $v$ is an endpoint of a path in ${\cal Q}$
immediately after $v$ is visited by the algorithm.


A pseudocode of the algorithm is given below (Algorithm~\ref{algo:SPC-tree}). In the description of the algorithm, we denote by $P_1+P_2$ the path obtained as the union of two edge-disjoint paths $P_1$ and $P_2$ meeting at a vertex. Similarly, $P_1+ P_2+ P_3$ denotes the union of three paths $(P_1+P_2)+P_3$.

{
\begin{algorithm}[h!]
\scriptsize
\caption{Solving the Shortest Path Collection problem in trees}
\label{algo:SPC-tree}
\KwIn{A tree $T= (V,E)$ and a non-empty set $A\subseteq V$ with $|A|$ even.}
\KwOut{A pair $({\cal A}, {\cal P})$, where ${\cal A}$ is a partition of the elements of $A$ into pairs $\{a_1,a_2\},\ldots, \{a_{2r-1},a_{2r}\}$,
minimizing $\sum_{i = 1}^r{\rm dist}_T(a_{2i-1},a_{i})$, and ${\cal P}$ is the collection of
corresponding shortest paths.}
\BlankLine
{
    Fix a root $r\in V$, and let $v_1,\ldots, v_n=r$ be the vertices of $T$ listed
    in reverse order with respect to the time they are visited by a breadth-first traversal from $r$.

    Set ${\cal A} = {\cal P} = {\cal Q} = \emptyset$.

    \For{$i= 1,2,\ldots, n$}
    {\nllabel{line-for}

        \If{$i<n$}
        {
            Set $q(v_i) = 0$.

            Let $C(v_i)$ be the set of children of $v_i$, and let $R(v_i) = \{u\in C(v_i)\,:\,q(u) = 1\}$.

            Fix an ordering $u_{1},\ldots, u_{k}$ of the elements of $R(v_i)$.


            \For{$j = 1,\ldots, \lfloor k/2\rfloor$}
            {
                Let $Q$ and $Q'$ be the paths in ${\cal Q}$ with $b(Q) = u_{2j-1}$ and $b(Q') = u_{2j}$.

                Add the path $Q+(u_{2j-1}, v_i, u_{2j})+Q'$ to ${\cal P}$.

                Add the pair $\{a(Q), a(Q')\}$ to ${\cal A}$.

                Remove $Q$ and $Q'$ from ${\cal Q}$.
            }

            \If{$k$ is odd}
            {
                Let $Q\in {\cal Q}$ be the path in ${\cal Q}$ with $b(Q) = u_k$.

                \If{$v_i\in A$}
                {
                    Add the path $Q+(u_k,v_i)$ to ${\cal P}$.

                    Add the pair $\{a(Q), v_i\}$ to ${\cal A}$.
                }
                \Else
                {
                    Add the path $Q'= Q + (u_k,v_i)$ to ${\cal Q}$, with $a(Q') = a(Q)$, $b(Q') = v_i$.

                    Set $q(v_i) = 1$.
                }
                Remove $Q$ from ${\cal Q}$.
            }
            \ElseIf{$v_i\in A$}
            {
                Add the trivial path $Q = (v_i)$ to ${\cal Q}$, with $a(Q) = b(Q) = v_i$.

                Set $q(v_i) = 1$.
            }
        }
        \Else{
            \tcp{we are at the root}
            Let $C$ be the set of children of $v_n = r$, and let $R  = \{v\in C\,:\,q(v) = 1\}$.

            Fix an ordering $u_{1},\ldots, u_{k}$ of the elements of $R$.

            \For{$j = 1,\ldots, \lfloor k/2\rfloor$}
            {
                Let $Q$ and $Q'$ be the paths in ${\cal Q}$ with $b(Q) = u_{2j-1}$ and $b(Q') = u_{2j}$.

                Add the path $Q+(u_{2j-1}, r, u_{2j})+Q'$ to ${\cal P}$.

                Add the pair $\{a(Q), a(Q')\}$ to ${\cal A}$.

                Remove $Q$ and $Q'$ from ${\cal Q}$.
            }
            \If{$k$ is odd}
            {
                \tcp{it must be the case that $r\in A$}
                Let $Q\in {\cal Q}$ be the path in ${\cal Q}$ with $b(Q) = u_k$.

                Add the path $Q +(u_k,r)$ to ${\cal P}$.

                Add the pair $\{a(Q), r\}$ to ${\cal A}$.

                Remove $Q$ from ${\cal Q}$.
            }
        }
    }
    \Return{(${\cal A}, {\cal P})$}
}
\end{algorithm}
}

To establish the correctness of the algorithm, we will show that the obtained solution attains a lower bound on the optimal value.
For every vertex $v$ of tree $T$ rooted at a fixed vertex $r$, let $k(v)$ denote the number of subtrees $T_i$ of $T$ rooted at the children of $v$ such that
$|A\cap V(T_i)|$ is odd. Then, for every feasible solution ${\cal A}'$ to the problem, the corresponding collection ${\cal P}'$ of shortest paths contains at least $k(v)$ edges connecting $v$ to its children. In particular, this implies that the optimal value of $\sum_{i = 1}^r{\rm dist}_T(a_{2i-1},a_{i})$
is at least $\sum_{v\in V(T)}k(v)$.

Now let us verify that the value of $\sum_{i = 1}^r{\rm dist}_T(a_{2i-1},a_{i})$ attained by the solution ${\cal A}$ given by Algorithm~\ref{algo:SPC-tree}
is equal to $\sum_{v\in V(T)}k(v)$. For every $i = 1,\ldots, n$,
let ${\cal P}_i$ and ${\cal Q}_i$ denote the collections of paths ${\cal P}$ and ${\cal Q}$
after $i$ iterations of the {\bf for} loop in line~\ref{line-for} (that is, immediately after vertex $v_i$
has been visited).
%
Properties $(1)$--$(4)$ described above (with ${\cal P}_i$ and ${\cal Q}_i$ in place of ${\cal P}$ and ${\cal Q}$, respectively)
can be proved by induction on $i$.
Moreover, for every $1\le i<j\le n$, it holds that ${\cal P}_i\subseteq {\cal P}_j$, every path in ${\cal Q}_i$ is a subpath of some path in ${\cal P}_j\cup {\cal Q}_j$, and
${\cal P}_n = {\cal P}$, where ${\cal P}$ is the final solution output by the algorithm.
Since the paths in ${\cal P}$ are edge-disjoint, the obtained value of $\sum_{i = 1}^r{\rm dist}_T(a_{2i-1},a_{i})$
is equal to the total number of edges that appear in paths in ${\cal P}$. Induction on $i$ can be used to show that:
\begin{enumerate}[(A)]
  \item After vertex $v_i$ has been processed, $q(v_i) = 1$ if and only if $|A\cap V(T_{v_i})|$ is odd, where $T_{v_i}$ is the subtree of $T$ rooted at $v_i$.
  \item The number of edges connecting $v_i$ to one of its children that are contained in a path from ${\cal P}_i\cup {\cal Q}_i$
  is equal to $k(v_i)$.
\item The number of edges connecting $v_i$ to one of its children that are contained in a path from ${\cal P}$
is the same as the number of edges connecting $v_i$ to one of its children that are contained in a path from ${\cal P}_i\cup {\cal Q}_i$.
\end{enumerate}
Since the final value of $\sum_{i = 1}^r{\rm dist}_T(a_{2i-1},a_{i})$ is equal to the total number of edges that appear in paths in ${\cal P}$, this value can be obtained by summing up, over all vertices $v_i$, the number
of edges connecting $v_i$ to its children that are contained in a path from ${\cal P}$.
By $(B)$ and $(C)$, this value is equal to the lower bound $\sum_{i = 1}^nk(v_i)$. Hence, the proof of correctness is complete.

It remains to analyze the algorithm's time complexity. A breadth-first traversal from $r$ takes $O(n)$ time.
A linear time implementation of the iterations of the {\bf for} loop can be achieved using an appropriate data structure representing the collections of
paths ${\cal P}$ and ${\cal Q}$ and their endpoints, and updating it at every iteration of the {\bf for} loop. (Each path in ${\cal P}\cup {\cal Q}$ can be represented by a doubly linked list.) The number of operations performed by the algorithm during  the $i$-th iteration of the {\bf for} loop  is then proportional to the degree of $v_i$. Altogether, this results in linear time complexity.
\end{proof}

The last thing to calculate is the number of operations needed to label the edges. We have to take into account three issues.

The first one is the division of $\gr$ into three sets: a one-element set with the identity element, the set of involutions and the set of pairs $\{a_i,-a_i\}$ of the remaining elements. In order to do that first we check the parity of all the cyclic subgroups $\gr_1,\ldots, \gr_k$, what can be done in time $O(k)=O(\log_2(|\gr|))$. Assume that the groups $\gr_1,\dots,\gr_p$ have even order and the groups $\gr_{p+1},\dots,\gr_k$ are odd for some $1\leq p\leq k$. Each element of $\gr$ is represented by some $k$-tuple $(g_1,\dots,g_k)$, where $g_j\in \gr_j$, $1\leq j\leq k$ and $0\leq g_j\leq |\gr_j|-1$. Then the identity element of $\gr$ is represented by $0_\gr=(0,\dots,0)$. The involutions have the form $(g_1,\dots,g_p,0,\dots,0)$, where $g_j\in\{0,|\gr_j|/2\}$ for $1\leq j\leq p$. The number of involutions is equal to $2^p-1$. Finally, the pairs $\{a_i,-a_i\}$ have the form $\{(g_1,\dots,g_p,g_{p+1},\dots,g_k),(|\gr_1|-g_1,\dots,|\gr_p|-g_p,|\gr_{p+1}|-g_{p+1},\dots,|\gr_k|-g_k)\}$, where $g_j\in\{0,\dots,|\gr_j|-1\}$ for $1\leq j\leq k-1$, $g_k\in\{0,\dots,\lfloor|\gr_k|\rfloor/2\}$ and for at least one $j$, $g_j\not\in I_j$, where $I_j=\{0,|\gr_j|/2\}$ for $1\leq j\leq p$ and $I_j=\{0\}$ for $p+1\leq j\leq k$. It is easy to observe that no search is necessary and the number of assignments is exactly $|\gr|$.

The second issue is to find (if necessary) the set of involutions (plus identity element in some cases) that sum up to $0$. The constructive proof of the Lemma \ref{lemma_inv} (see \cite{ref_AnhCic2}) gives the simple method of choosing such a subset, with time complexity linear in the number of involutions. If the required number of elements $r\in\{0,1,2^p-1,2^p\}$ then $R=\emptyset$ or $R=\{0\}$ or $R=I\setminus \{0\}$ or $R=I$, respectively. If $3\leq r\leq 2^{p-1}$, then we select $r-1$ elements one by one in the lexicographic order (here $0_j$ denotes the identity element of $\gr_j$, in other words the $0$ on $j_th$ position of the $k$-tuple): $(0_1,\dots,0_p,0_{p+1},\dots,0_k)$, $(0_1,\dots,|\gr_p|/2,0_{p+1},\dots,0_k)$, $(0_1,\dots,|\gr_{p-1}|/2,0_p,0_{p+1},\dots,0_k)$, $(0_1,\dots,|\gr_{p-1}|/2,|\gr_{p}|/2,0_p,0_{p+1},\dots,0_k)$ and so on. Now we calculate the sum of all chosen elements. It has the form $s=(0_1,g_2\dots,g_p,0_{p+1},\dots,0)$, similarly to all the elements of the list. If $s$ is not on the list, then we add it to the list and we are done. Otherwise we choose another element $s_1$ of the list, (if $s\neq 0_\gr$ then the one preceding $s$, otherwise the one following it), and we change the first coordinate of both $s$ and $s_1$ from $0_1$ to $|\gr_1|/2$ and we are done. The last case is when $2^{p-1}+1\leq r\leq 2^p-3$. In such a case we construct the $(2^p-r)$-element set $R_1$ and then take $R=I\setminus R_1$. As it can be easily seen, the complexity is at most $r-1$ (choosing the elements) plus $r-2$ (sum) plus $r-1$ (checking if $s$ is on the list) plus eventually $2$ additional summations, what gives in total $O(|\gr|)$.

The third issue is the assignment of labels to the paths. If $G$ is not a star, then it is equal to the sum of the numbers of edges in the optimal solutions of the {\sc Shortest Path Collection} for $A=V_0$, $A=V_1\setminus V_0$ and $A=V_2\setminus V_0$, where $V_0$ is either an empty set or some subset of $V(G)$ with at most $8$ elements, while $V_1$ and $V_2$ are the color classes of the spanning tree $T(G)$ of $G$. As all the paths in the optimal solution of {\sc Shortest Path Collection} are edge disjoint and we assign $0$ to all the edges in $E(G)\setminus E(T(G))$, the total number of assignments does not exceed $2|E(T(G))|+|E(G)|< 2|V(G)|+|E(G)|$.

All the above calculations make sense if $|\gr|$ fulfills all the necessary conditions. Obviously they can be checked in constant time. Hence the following corollary is true.

\begin{mycorollary}
The {\sc Irregular Labeling}  problem can be solved in time $O(|E(G)|+|\gr|)$.
\end{mycorollary}

Observe that in the case of $|\gr|=s_g(G)$ the complexity reduces to $O(|E(G)|)$. Note also that instead of the list of orders of prime-power cyclic groups, $\gr$ can be represented with the minimum-length list of the orders of cyclic groups (e.g. $(2,6)$ instead of $(2,2,3)$, as $\mathbb{Z}_6\cong \mathbb{Z}_2\times \mathbb{Z}_3$), and it would not change the method of the division of $\gr$. It can also be represented as the list of generators with the relations of the form $n_jg_j=0$ (in the above example: $\gr=[g_1,g_2|2g_1=6g_2=0]$) but in such a case the list of multipliers $n_j$ is equivalent to the list of orders of cyclic groups. Of course, in all those cases the complexity of division remains $O(|\gr|)$.

\section{Final Remarks}

As we solved the problem for all connected graphs, next step would be to find the solution for arbitrary graphs.

\begin{myproblem}
Determine  group irregularity strength $s_g(G)$ for disconnected graph $G$ with no component of order less than $3$.
\end{myproblem}

Theorem~\ref{main_thm2} characterises all the pairs $(G,\gr)$, $|\gr|>s_g(G)$, such that there exists an irregular-$\gr$-labeling of given connected graph $G$ of order at least $3$. Thus the following generalisation arises.

\begin{myproblem}
Characterise all the pairs $(G,\gr)$, $|\gr|>s_g(G)$, such that there exists an irregular-$\gr$-labeling of given graph $G$ with no components of order less than $3$.
\end{myproblem}

In the proof of Theorem~\ref{main_thm} we often use the fact that we are allowed to use $0$ on edges. Thus next natural problem is the following.

\begin{myproblem}
Let $G$ be a simple graph with no components of order less than $3$. For any Abelian group $\gr$, let $\gr^* = \gr \setminus \{0\}$. Determine \textit{non-zero group irregularity strength} ($s^*_g(G)$) of $G$, i.e., the smallest value of $s$ such that taking any Abelian group $\gr$ of order $s$, there exists a function $f:E(G)\rightarrow \gr^*$ such that the sums of edge labels in every vertex are distinct.
\end{myproblem}

All the elements of $\gr$ can be obtained as some combination of not necessarily all of its elements, in particular of its generators. The question is, how many elements of $\gr$ we have to use in order to obtain $\gr$-irregular labelling.

\begin{myproblem}
Assume that for given simple graph $G$ with no components of order less than $3$ there exists $\gr$-irregular labelling for every group $\gr$ of order $s$. What is the minimum number $k=k(G,s)$ such that for every group $\gr$ of order $s$ there is a subset $S\subseteq\gr$, $|S|\leq k$ such that there exists a $\gr$-irregular labelling $f:E(G)\rightarrow S$?
\end{myproblem}

So far we considered only finite Abelian groups. So, next question seems to be natural, as some generalization of the problem of the ordinary irregularity strength.

\begin{myproblem}
Let $G$ be a simple graph with no component of order less than $3$. Determine the smallest value of $k$ such that for any infinite Abelian group $\gr$ there exists a subset $S\subseteq\gr$, $S\leq k$ such that there exists a $\gr$-irregular labelling $f:E(G)\rightarrow S$.
\end{myproblem}

\subsection*{Acknowledgments}
M.~Anholcer and S.~Cichacz are grateful to the Leaders, Researchers and Administrative Staff of the Faculty of Mathematics, Natural Sciences and Information Technologies of the University of Primorska for creating the perfect work conditions during their post-doc fellowship and inspiring discussions. We are also grateful to Prof.~Dalibor Froncek for inspiration. M.~Milani\v c is supported in part by ``Agencija za raziskovalno dejavnost Republike Slovenije'', research program P$1$--$0285$ and research projects J$1$--$4010$, J$1$--$4021$ and N$1$--$0011$.


\begin{thebibliography}{99}

\bibitem{ref_AigTri}
Aigner M., Triesch E., \textit{Irregular assignments of trees and forests}, SIAM Journal on Discrete Mathematics Vol.3 No.4 (1990), 439-449.

\bibitem{ref_AmaTog}
Amar D., Togni O., \textit{Irregularity strength of trees}, Discrete Mathematics 190 (1998), 15-38.

\bibitem{ref_AnhCic2}
Anholcer M., Cichacz S., \emph{Group sum chromatic number of graphs}, preprint, arXiv:1205.2572v1 [math.CO].

\bibitem{ref_BeaGalHeaJun}
Beals R., Gallian J., Headley P., Jungreis D., \textit{Harmonious groups}, Journal of Combinatorial Theory, Series A, 56 (1991) pp. 223-238.

\bibitem{ref_CavComNel}
Cavenagh N., Combe D., Nelson A. M., \textit{Edge-magic group labellings of countable graphs}, Electronic Journal of Combinatorics, 13 (2006) {\#}R92, 19 pp.

\bibitem{ref_ChaJacLehOelRuiSab1}
Chartrand G., Jacobson M.S., Lehel J., Oellermann O.R., Ruiz S., Saba F., \textit{Irregular networks}, Congressus Numerantium 64 (1988), 187-192.

\bibitem{ref_ComNelPal}
Combe D., Nelson A.M., Palmer W.D., \textit{Magic labellings of graphs over finite Abelian groups}, Australasian Journal of Combinatorics 29 (2004), pp. 259-271.

\bibitem{ref_Edm}
Edmonds J., \textit{Maximum Matching and a Polyhedron with 0,1-Vertices}, Journal of Research of the National Bureau of Standards 69B (1965) 125--130.

\bibitem{ref_Fro}
Froncek D., \textit{Group distance magic labeling of $C_k \square C_m$}, preprint.

\bibitem{ref_FerGouKarPfe}
Ferrara M., Gould R.J., Karo\'nski M., Pfender F.: \textit{An iterative approach to graph irregularity strength}, Discrete Applied Mathematics 158 (2010), pp. 1189-1194.

\bibitem{ref_GraSlo}
Graham R.L., Sloane N. J. A., \textit{On additive bases and harmonious graphs}, SIAM Journal
on Algebraic and Discrete Methods, 1 (1980), pp. 382-404.

\bibitem{ref_Hov}
Hovey M., \textit{$A$-cordial graphs}, Discrete Mathematics 93 (1991), pp. 183-194.

\bibitem{ref_KalKarPfe1}
Kalkowski M., Karo\'nski M., Pfender F., \textit{A new upper bound {for~the~irregularity} strength of graphs}, SIAM Journal on Discrete Mathematics 25 (2011), pp. 1319-1321.

\bibitem{ref_KapLevRod}
Kaplan G., Lev A., and Roditty Y., \textit{On zero-sum partitions and anti-magic trees}, Discrete
Mathematics 309 (2009), pp. 2010-2014.

\bibitem{ref_Mih}
Mih\v{a}ilescu P., \textit{Primary cyclotomic units and a proof of Catalan's Conjecture}, Journal f\"{u}r die reine und angewandte Mathematik 572 (2004), pp. 167-195.

\bibitem{ref_Prz3}
Przyby{\l}o J., personal communication.

\bibitem{ref_Sta}
Stanley R. P., \textit{Linear homogeneous Diophantine equations and magic labelings of graphs},
Duke Mathematical Journal 40 (1973), pp. 607-632.

\bibitem{ref_Tog1}
Togni O., \textit{Force des graphes. Indice optique des r\'{e}seaux}, Th\`{e}se pr\'{e}sent\'{e}e pour obtenir le grade de docteur, Universit\'{e} de Bordeaux 1, \'{E}cole doctorale de mathematiques et d'informatique (1998), 141 pp.

\bibitem{ref_Zak}
\.Zak A., \textit{Harmonious orders of graphs}, Discrete Mathematics 309 (2009), pp.6055-6064.

\end{thebibliography}
\end {document}